\documentclass[a4paper,12pt]{article}
\usepackage{amssymb,amscd,amsfonts,amsbsy}
\usepackage{enumerate}
\usepackage{epsfig}
\input amssymb.sty
\input euscript.sty

\def\ring{\mathaccent"0017 }

\newcommand{\RR}{{\mathbb{R}}}
\newcommand{\NN}{{\mathbb{N}}}

\newcommand{\po}{\partial\Omega}
\def\ring{\mathaccent"0017 }

\newcommand{\bp}{\noindent {\it Proof}.\,\,}
\newcommand{\ep}{\hfill$\Box$ \vskip 0.08in}

\newcommand{\meanint}{{\int{\mkern-19mu}-}}

\newtheorem{proposition}{Proposition}[section]
\newtheorem{theorem}[proposition]{Theorem}
\newtheorem{lemma}[proposition]{Lemma}
\newtheorem{corollary}[proposition]{Corollary}

\begin{document}

\title{Boundedness of the Hessian of a biharmonic
function in a convex domain}

\author{Svitlana Mayboroda and Vladimir Maz'ya\thanks{
2000 {\it Math Subject Classification:} 35J40, 35J30, 35B65.
\newline The second author is partially supported by NSF grant DMS 0500029.}}

\date{ }

\maketitle

\begin{abstract} We consider the
Dirichlet problem for the biharmonic equation on an arbitrary
convex domain and prove that the second derivatives of the
variational solution are bounded in all dimensions.
\end{abstract}

\section{Introduction}
\setcounter{equation}{0}

Properties of solutions of the Dirichlet problem for the
Laplacian on convex domains are nowadays well understood. It is a
classical fact  that the gradient of a solution is bounded  and in
the last decade a number of results in $L^p$, Sobolev  and Hardy
spaces have been developed (see  \cite{AdL2}, {\cite{AdLp},
\cite{MaMi}, \cite{Fr1}, {\cite{FrJe}). However, much less is
known about the behavior of  solutions to higher order elliptic
equations. The aim of this paper is to establish the boundedness
of the second derivatives of a  biharmonic function in all
dimensions.

To be more precise, given a bounded domain $\Omega\subset\RR^n$
denote by $\ring W_2^2(\Omega)$ the completion of
$C_0^\infty(\Omega)$ in the norm of the Sobolev space of functions
with second distributional derivatives in $L^2$. We consider the
variational solution of the boundary value problem

\begin{equation}\label{eq1.1}
\Delta^2 u=f \,\,{\mbox{in}}\,\,\Omega, \quad f\in
C_0^{\infty}(\Omega),\quad u\in \ring W_2^2(\Omega).
\end{equation}

\noindent The main result of this paper is the following.

\begin{theorem}\label{t1.1}
Let $\Omega$ be a convex domain in $\RR^n$, $O\in\po$, and fix
some $R\in(0,{\rm diam}\,(\Omega)/10)$. Suppose $u$ is a solution
of the Dirichlet problem {\rm{(\ref{eq1.1})}} with $f\in
C_0^{\infty}(\Omega\setminus B_{10R})$. Then

\begin{equation}\label{eq1.3}
|\nabla^2 u(x)| \leq
\frac{C}{R^2}\,\left(\meanint_{C_{R/2,\,5R}\cap\Omega}
|u(x)|^2\,dx\right)^{1/2}\quad \mbox{for every }\quad x\in
B_{R/5}\cap\Omega,
\end{equation}

\noindent where $\nabla^2u$ is the Hessian matrix of $u$,

\begin{equation}\label{eq1.4}
C_{\rho,\,R}=\{x\in\RR^n:\,\rho\leq|x|\leq R\},\quad
B_{\rho}=\{x\in\RR^n:\,|x|< \rho\},
\end{equation}

\noindent and the constant $C$ depends on the dimension only.

In particular,

\begin{equation}\label{eq1.4.1}
|\nabla^2 u| \in L^\infty(\Omega).
\end{equation}
\end{theorem}

We would like to mention that the properties of solutions to
boundary value problems for the biharmonic equation on general
domains, such as pointwise estimates and an analogue of the Wiener
criterion, have been studied in  \cite{M1}, \cite{M11}, \cite{M2}. 
The boundedness of the Hessian of a
biharmonic function on a planar
convex domain  was obtained from the asymptotic
formulas in \cite{KM}, where the restriction $n=2$ was essential.
In the
context of Lipschitz domains, following the work of B.\,Dahlberg,
C.\,Kenig and G.\,Verchota in \cite{DKV} on the well-posedness of
the Dirichlet problem with boundary data in $L^2$, J.\, Pipher and
G.\, Verchota established
 the so-called boundary G\aa
rding inequality (\cite{PVGarding}), $L^p$ estimates (\cite{PVLp})
and the Miranda-Agmon maximum principle in low dimensions
(\cite{PVmax}). The more recent advances include the work of Z.\,
Shen (\cite{Shen0}, \cite{Shen}) and regularity results by
V.\,Adolfsson and J.\, Pipher in (\cite{PAd}).

The
result of Theorem~\ref{t1.1} may fail on a Lipschitz domain. Even
in the three-dimensional case the solution in the exterior of a
thin cone is only $C^{1,\alpha}$ for some $\alpha>0$ (see
\cite{KMR}) and there is a four-dimensional domain for which the
gradient of the solution is not bounded (see \cite{MR} and
\cite{PVLp} for counterexamples).

Our approach is different from the methods of the Lipshitz theory
as well as from those used in \cite{KM}. The proof of (\ref{eq1.4.1}) relies 
upon the new weighted integral identities, which allow  to control the
local $L^2$ behavior of biharmonic functions near the boundary of
the domain.

\section{Global estimates: part I}
\setcounter{equation}{0}

Let $\Omega$ be an arbitrary domain in $\RR^n$, $n\geq 2$. We
assume that the origin belongs to the complement of $\Omega$ and
$r=|x|$, $\omega=x/|x|$ are the spherical coordinates centered at
the origin. In fact, we will mostly use the coordinate system
$(t,\omega)$, where $t=\log r^{-1}$, and the mapping $\varkappa$
defined by

\begin{equation}\label{eq3.1}
\RR^n \ni x\,\stackrel{\varkappa}{\longrightarrow}
\,(t,\omega)\in\RR\times S^{n-1}.
\end{equation}

\noindent  Here, and throughout the paper, $S^{n-1}$ denotes the
unit sphere in $\RR^n$,

\begin{equation}\label{eq3.1.1}
\RR^n_+=\{x\in\RR^n:\,x_n>0\},
\end{equation}

\noindent and $S^{n-1}_{+}=S^{n-1}\cap \RR^n_+$.

Next, given an open set $\Gamma\subset S^{n-1}$, the space $\ring
W_2^1(\Gamma)$ is the completion of $C_0^\infty(\Gamma)$ in the
norm

\begin{equation}\label{eq2.8}
\|\psi\|_{\ring W_2^1(\Gamma)}=\left(-\int_{\Gamma}\bar
\psi\,\delta_\omega
\psi\,d\omega\right)^{1/2}=\left(\int_{\Gamma}|\nabla_\omega
\psi|^2\,d\omega\right)^{1/2},
\end{equation}

\noindent where $\delta_\omega$ and $\nabla_\omega$ stand for the
Laplace-Beltrami operator and gradient on the unit sphere,
respectively, and $\ring W_2^2(\Gamma)$ is the completion of
$C_0^\infty(\Gamma)$ with respect to the norm

\begin{equation}\label{eq2.19}
\|\psi\|_{\ring W_2^2(\Gamma)}=\left(\int_\Gamma|\delta_\omega
\psi|^2\,d\omega\right)^{1/2}.
\end{equation}

If $\Gamma\subset S^{n-1}_+$ and $\psi\in \ring W_2^k(\Gamma)$,
$k=1,2$, we sometimes write that $\psi$ belongs to $\ring
W_2^k(S^{n-1}_+)$, with the understanding that $\psi$ is extended
by zero outside of $\Gamma$, and similarly the functions
originally defined on the subsets of $\RR^n$ will be extended by
zero and treated as functions on $\RR^n$ whenever appropriate.

\begin{lemma}\label{l3.1}
Let $\Omega$ be an arbitrary bounded domain in $\RR^n$,
$O\in\RR^n\setminus\Omega$ and

\begin{equation}\label{eq3.2}
u\in C^2(\bar\Omega),\quad u\Big|_{\po}=0, \quad \nabla
u\Big|_{\po}=0, \quad v=e^{2 t}(u\circ \varkappa^{-1}).
\end{equation}

\noindent  Then

\begin{eqnarray}\label{eq3.3}
&&\hskip -1cm \int_{\Omega}\Delta
u(x)\Delta\left(\frac{u(x)G(\log |x|^{-1})}{|x|^{n}}\right)\,dx\nonumber\\[4pt]
&&\hskip -1cm\quad\quad  =\int_{\varkappa(\Omega)}\Bigl[
(\delta_\omega v)^2G+2(\partial_t\nabla_\omega v)^2G+
(\partial_t^2v)^2G\nonumber\\[4pt]
&&\hskip -1cm\quad\quad    -(\nabla_\omega v)^2 \Bigl(\partial_t^2
G+n\,\partial_tG+2n\,G\Bigr) -(\partial_t
v)^2\Bigl(2\partial_t^2G+3n\,\partial_t
G+\bigl(n^2+2n-4\bigr)\,G\Bigr)
\nonumber\\[4pt]
&&\hskip -1cm\quad\quad   +\frac
12\,v^2\Bigl(\partial_t^4G+2n\,\partial_t^3G
+\bigl(n^2+2n-4\bigr)\,\partial_t^2G+
2n(n-2)\,\partial_tG\Bigr)\Bigr]\,d\omega dt.
\end{eqnarray}

\noindent for every  function $G$ on $\RR$ such that both sides of
{\rm{(\ref{eq3.3})}} are well-defined.
\end{lemma}

\vskip 0.08in

\bp In the coordinates $(t,\omega)$ the $n$-dimensional Laplacian
can be represented as

\begin{equation}\label{eq2.5}
\Delta=e^{2t}\Lambda(\partial_t,\delta_\omega),
\quad\mbox{where}\quad
\Lambda(\partial_t,\delta_\omega)=\partial_t^2-(n-2)\partial_t+\delta_\omega.
\end{equation}

\noindent  Then

\begin{eqnarray}\label{eq3.4}
&&\int_{\Omega}\Delta
u(x)\Delta\left(\frac{u(x)G(\log |x|^{-1})}{|x|^{n}}\right)\,dx\nonumber\\[4pt]
&&\quad = \int_{\varkappa(\Omega)}
\Lambda(\partial_t-2,\delta_\omega)
v\,\Lambda(\partial_t+n-2,\delta_\omega)(vG)\,d\omega dt\nonumber\\[4pt]
&&\quad = \int_{\varkappa(\Omega)}
\left(\partial_t^2v-(n+2)\partial_tv+2nv+\delta_\omega v\right)
\,\left(\partial_t^2(vG)+(n-2)\partial_t(vG)+G\,\delta_\omega
v\right)\,d\omega dt\nonumber\\[4pt]
&&\quad = \int_{\varkappa(\Omega)}
\left(\partial_t^2v-(n+2)\partial_tv+2nv+\delta_\omega v\right)\nonumber\\[4pt]
&&\quad\quad \times\left(G\,\delta_\omega
v+G\,\partial_t^2v+(2\partial_t
G+(n-2)G)\,\partial_tv+(\partial_t^2G+(n-2)\partial_t
G)\,v\right)\,d\omega dt.
\end{eqnarray}

\noindent Expanding the expression above and reassembling the
terms, we write the integral in (\ref{eq3.4}) as

\begin{eqnarray}\label{eq3.5}
&& \int_{\varkappa(\Omega)} \Bigl( G\,(\delta_\omega
v)^2+2G\,\delta_\omega
v\partial_t^2 v+G\,(\partial_t^2v)^2\nonumber\\[4pt]
&&\quad\quad+v\delta_\omega v\,\left(\partial_t^2
G+(n-2)\partial_tG+2nG\right)+\delta_\omega v\partial_t
v\,\left(2\partial_t G-4G\right)
\nonumber\\[4pt]
&&\quad\quad +\partial_t^2 v\partial_t
v\left(2\partial_tG-4G\right)+v\partial_t^2
v\left(\partial_t^2G+(n-2)\partial_tG+2nG\right)\nonumber\\[4pt]
&&\quad\quad+(\partial_tv)^2\,\left(-2(n+2)\partial_t
G-(n^2-4)G\right)
\nonumber\\[4pt]
&&\quad\quad+v\partial_t
v\left(-(n+2)\partial_t^2G-(n^2-4n-4)\partial_tG+2n(n-2)G\right)\nonumber\\[4pt]
&&\quad\quad+v^2\left(2n\partial_t^2G+2n(n-2)\partial_tG\right)
\Bigr)\,d\omega dt.
\end{eqnarray}

\noindent Since $G$ does not depend on $\omega$, after integration
by parts the latter integral becomes

\begin{eqnarray}\label{eq3.6}
&&\hskip -0.7cm  \int_{\varkappa(\Omega)} \Bigl( G\,(\delta_\omega
v)^2-2G\,\delta_\omega\partial_t
v\partial_t v+G\,(\partial_t^2v)^2\nonumber\\[4pt]
&&\hskip -0.7cm\quad+(\nabla_\omega
v)^2\,\left(-\partial_t^2G-\partial_t^2
G-(n-2)\partial_tG-2nG+\partial_t^2 G-2\partial_tG\right)
\nonumber\\[4pt]
&&\hskip -0.7cm\quad +(\partial_t
v)^2\left(-\partial_t^2G+2\partial_tG-\left(\partial_t^2G+(n-2)\partial_tG+2nG\right)-2(n+2)\partial_t
G-(n^2-4)G\right)
\nonumber\\[4pt]
&&\hskip -0.7cm\quad+v\partial_t
v\left(-\partial_t^3G-(n-2)\partial^2_tG-2n\partial_tG
-(n+2)\partial_t^2G-(n^2-4n-4)\partial_tG+2n(n-2)G\right)\nonumber\\[4pt]
&&\hskip
-0.7cm\quad+v^2\left(2n\partial_t^2G+2n(n-2)\partial_tG\right)
\Bigr)\,d\omega dt.
\end{eqnarray}

\noindent Finally, integrating by parts once again and collecting
the terms, we arrive at (\ref{eq3.3}). \ep

To proceed further we need some auxiliary results. By $\delta$ we
denote the Dirac delta function.

\begin{lemma}\label{l3.2} A bounded solution of the equation

\begin{equation}\nonumber
\frac{d^4g}{dt^4}+2n\frac{d^3g}{dt^3}
+\left(n^2-2\right)\frac{d^2g}{dt^2}-2n\frac{dg}{dt} =\delta
\label{eq3.7}
\end{equation}

\noindent subject to the restriction

\begin{equation}\label{eq3.8}
g(t)\to 0 \mbox{ as } t \to +\infty,
\end{equation}

\noindent is the function

\begin{equation}\label{eq3.9}
g(t)=-\frac 1{2n\sqrt{n^2+8}}\left\{\begin{array}{l}
n\, e^{- 1/2(n-\sqrt{n^2+8})t}-\sqrt{n^2+8},\qquad \qquad\qquad t<0,\\[4pt]
n\, e^{ -1/2(n+\sqrt{n^2+8})t} -\sqrt{n^2+8}\,e^{-nt},\qquad\qquad t>0.\\[4pt]
\end{array}
\right.
\end{equation}

\end{lemma}

\bp The equation (\ref{eq3.7}) can be written as

\begin{equation}\label{eq3.10}
\frac{d}{dt}\left(\frac{d}{dt}+n\right)\left(\frac{d}{dt}+\frac
12\left(n+\sqrt{n^2+8}\right)\right) \left(\frac{d}{dt}+\frac
12\left(n-\sqrt{n^2+8}\right)\right)g=\delta.
\end{equation}

\noindent Since we seek a bounded solution of (\ref{eq3.7})
satisfying (\ref{eq3.8}), $g$ must have the form

\begin{equation}\label{eq3.9.1}
g(t)=\left\{\begin{array}{l}
a\, e^{- 1/2(n-\sqrt{n^2+8})t}+b,\qquad \qquad\qquad t<0,\\[4pt]
c\, e^{ -1/2(n+\sqrt{n^2+8})t} +d\,e^{-nt},\qquad\qquad t>0,\\[4pt]
\end{array}
\right.
\end{equation}

\noindent for some constants $a,b,c,d$. Once this is established,
we find the system of coefficients so that $\partial_t^kg$ is
continuous for $k=0,1,2$ and $\lim_{t\to 0^+}\partial_t^3
g(t)-\lim_{t\to 0^-}\partial_t^3 g(t)=1$.
 \ep

Next, let us consider some estimates based on the spectral
properties of the Laplace-Beltrami operator on the half-sphere. When $n\geq 3$ we
will use the coordinates $\omega=(\theta,\varphi)$ on the unit
sphere $S^{n-1}$, where $\theta\in[0,\pi]$, $\varphi\in S^{n-2}$. In the two-dimensional case
$\omega=\theta\in [0,2\pi)$.

\begin{lemma}\label{l3.2.2} For every $v\in \ring W_2^2(S^{n-1}_+)$,

\begin{equation}\label{eq3.10.2}
\int_{S^{n-1}_+}|\delta_\omega v|^2\,d\omega\geq 2n
\int_{S^{n-1}_+}|\nabla_\omega v|^2\,d\omega.
\end{equation}

\noindent The equality is achieved when $v=\cos^2
\theta$.
\end{lemma}

\bp Since $v \in \ring W_2^2(S^{n-1}_+)$,

\begin{equation}\label{eq3.10.3}
\int_{S^{n-1}_+}|\nabla_\omega v|^2\,d\omega=
\int_{S^{n-1}_+}\left|\nabla_\omega(v-(v)_{S^{n-1}_+})\right|^2\,d\omega\leq
\|v-(v)_{S^{n-1}_+}\|_{L^2(S^{n-1}_+)} \|\delta_\omega
v\|_{L^2(S^{n-1}_+)},
\end{equation}

\noindent where

\begin{equation}\label{eq3.10.3.1}
(v)_{S^{n-1}_+}=\meanint_{S^{n-1}_+}v\,d\omega.
\end{equation}

\noindent If $n\geq 3$, let us denote

\begin{equation}\label{eq3.10.4}
z(\theta):=\meanint_{S^{n-2}}v(\theta,\varphi)\,d\varphi,\qquad
y(\theta,\varphi):=v(\theta,\varphi)-z(\theta).
\end{equation}

\noindent Then $(y)_{S^{n-1}_+}=0$ and

\begin{eqnarray}\label{eq3.10.5}
 \|v-(v)_{S^{n-1}_+}\|_{L^2(S^{n-1}_+)}^2
&=&\|y\|_{L^2(S^{n-1}_+)}^2+\|z-(z)_{S^{n-1}_+}\|_{L^2(S^{n-1}_+)}^2,
\\[4pt] \label{eq3.10.6}
 \|\nabla_\omega v\|_{L^2(S^{n-1}_+)}^2&=& \|\nabla_\omega
y\|_{L^2(S^{n-1}_+)}^2+\|\nabla_\omega z\|_{L^2(S^{n-1}_+)}^2.
\end{eqnarray}

By the definition (\ref{eq3.10.4}) the function $y\in \ring
W_2^1(S^{n-1}_+)$ is orthogonal to $\cos \theta$ on $S^{n-1}_+$.
Therefore, $y$ is orthogonal to the first eigenfunction of the
Dirichlet problem for $-\delta_\omega$ on $S^{n-1}_+$. Since the
second eigenvalue of $-\delta_\omega$ is $2n$, this yields

\begin{equation}\label{eq3.10.7}
\int_{S^{n-1}_+}|\nabla_\omega y|^2\,d\omega\geq 2n
\int_{S^{n-1}_+}|y|^2\,d\omega.
\end{equation}

Turning to the estimates on $z$, we observe that

\begin{equation}\label{eq3.10.8}
\Lambda=\inf\left\{\frac{-\int_{S^{n-1}_+}\bar \xi\,\delta_\omega
\xi\,d\omega}{\int_{S^{n-1}_+}|\xi|^2\,d\omega}: \quad \xi\in
W^1_2({S^{n-1}_+}),\,\,\int_{S^{n-1}_+}\xi\,d\omega=0,\,\,\xi=\xi(\theta)\right\},
\end{equation}

\noindent is the first positive eigenvalue of the Neumann problem
for $-\delta_\omega$ in the space of axisymmetric functions.
Hence, $\Lambda=2n$ and the corresponding eigenfunction is
$n\cos^2\theta-1$. Therefore,

\begin{equation}\label{eq3.10.9}
\int_{S^{n-1}_+}|\nabla_\omega z|^2\,d\omega\geq 2n
\int_{S^{n-1}_+}|z-(z)_{S^{n-1}_+}|^2\,d\omega.
\end{equation}

\noindent Combined with (\ref{eq3.10.5})--(\ref{eq3.10.7}),
the formula above implies that

\begin{equation}\label{eq3.10.9.0}
\int_{S^{n-1}_+}|\nabla_\omega v|^2\,d\omega\geq 2n
\int_{S^{n-1}_+}|v-(v)_{S^{n-1}_+}|^2\,d\omega,
\end{equation}

\noindent and by (\ref{eq3.10.3}) this
finishes the proof for $n\geq 3$.

In the case $n=2$ there is no need to introduce functions $z$ and $y$.
One can work directly with $v$ and prove (\ref{eq3.10.9.0})
following the argument for (\ref{eq3.10.9}) above.
\ep

\begin{lemma}\label{l3.3}
Let $\Omega$ be a bounded convex domain in $\RR^n$ and $O\in
\RR^n\setminus\Omega$. Suppose that

\begin{equation}\label{eq3.11}
u\in C^2(\bar\Omega),\quad u\Big|_{\po}=0, \quad \nabla
u\Big|_{\po}=0, \quad v=e^{2 t}(u\circ \varkappa^{-1}),
\end{equation}

\noindent  and $g$ is given by {\rm{(\ref{eq3.9})}}. Then

\begin{eqnarray}\label{eq3.12}
&&\hskip -1.5 cm \int_{\Omega}\Delta
u(x)\Delta\left(\frac{u(x)g(\log (|\xi|/|x|))}{|x|^{n}}\right)\,dx\nonumber\\[4pt]
&&\hskip -1.5 cm \quad
\geq-\int_{\varkappa(\Omega)}\Bigl(2\partial_t^2g(t-\tau)+3n\partial_t
g(t-\tau)+\left(n^2-2\right)g(t-\tau)\Bigr)(\partial_t
v(t,\omega))^2\,d\omega dt
\nonumber\\[4pt]
&&\hskip -1.5 cm \qquad+\frac 12\int_{S^{n-1}\cap
\varkappa(\Omega)}v^2(\tau,\omega)\,d\omega,
\end{eqnarray}

\noindent for every $\xi\in\Omega$, $\tau=\log |\xi|^{-1}$.

\end{lemma}

\bp Rearranging the terms in (\ref{eq3.3}), we write

\begin{eqnarray}\label{eq3.13}
&&\hskip -0.7cm \int_{\Omega}\Delta
u(x)\Delta\left(\frac{u(x)G(\log |x|^{-1})}{|x|^{n}}\right)\,dx\nonumber\\[4pt]
&&\hskip -0.7cm \quad =\int_{\varkappa(\Omega)}\Bigl[G \Bigl(
(\delta_\omega v)^2 -2n(\nabla_\omega v)^2 + (\partial_t^2v)^2
+2(\partial_t\nabla_\omega v)^2 -\bigl(n^2+2n-4\bigr)(\partial_t
v)^2\Bigr)\nonumber\\[4pt]
&&\hskip -0.7cm \qquad  -\Bigl(2\partial_t^2G+3n\partial_t
G\Bigr)(\partial_t v)^2-\Bigl(\partial_t^2
G+n\partial_tG\Bigr)(\nabla_\omega v)^2
\nonumber\\[4pt]
&&\hskip -0.7cm \qquad +\frac
12\Bigl(\partial_t^4G+2n\partial_t^3G
+\bigl(n^2+2n-4\bigr)\partial_t^2G
+2n(n-2)\partial_tG\Bigr)v^2\Bigr]\,d\omega dt.
\end{eqnarray}

\noindent Let $G(t)=g(t-\tau)$, $t\in\RR$, with $g$ defined in
(\ref{eq3.9}). First of all observe that for such a choice of $G$
the conclusion of Lemma~\ref{l3.1} (and hence the equality
(\ref{eq3.13})) remains valid.

Going further, $g\geq 0$ and

\begin{equation}\label{eq3.14}
\partial_t^2 g(t)+n\partial_tg(t)\leq 0,
\quad{\mbox{for every}}\quad t\in\RR.
\end{equation}

\noindent Indeed,

\begin{equation}\label{eq3.15}
\partial_t\,g(t)=\frac 1{2\sqrt{n^2+8}}\left\{\begin{array}{l}
 \frac{1}{2}\left(n-\sqrt{n^2+8}\right) e^{- 1/2(n-\sqrt{n^2+8})t},\qquad
 \qquad\qquad \quad\,\,\,\quad t<0,\\[4pt]
 -\sqrt{n^2+8}\, e^{-nt}+
\frac{1}{2}\left(n+\sqrt{n^2+8}\right)e^{-
1/2(n+\sqrt{n^2+8})t},\quad t>0,
\end{array}
\right.
\end{equation}

\noindent and

\begin{equation}\label{eq3.16}
\partial_t^2\,g(t)=\frac 1{2\sqrt{n^2+8}}\left\{\begin{array}{l}
- \frac{1}{4}\left(n-\sqrt{n^2+8}\right)^2 e^{- 1/2(n-\sqrt{n^2+8})t},\quad
\,\,\,\,\,\quad \qquad\qquad t<0,\\[4pt]
n \sqrt{n^2+8}\, e^{-nt}-
\frac{1}{4}\left(n+\sqrt{n^2+8}\right)^2e^{-
1/2(n+\sqrt{n^2+8})t},\,\, t>0.
\end{array}
\right.
\end{equation}

\noindent Therefore,

\begin{equation}\label{eq3.17}
\partial_t^2
g(t)+n\partial_tg(t)=-\frac{1}{\sqrt{n^2+8}}\left\{\begin{array}{l}
e^{- 1/2(n-\sqrt{n^2+8})t},\qquad t<0,\\[4pt]
e^{- 1/2(n+\sqrt{n^2+8})t},\qquad t>0,\\[4pt]
\end{array}
\right.
\end{equation}

\noindent is non-positive.

Recall that the first eigenvalue of the operator $-\delta_\omega$
on the half-sphere is $n-1$. Since $\Omega$ is a convex domain,
for every fixed $t\in\RR$

\begin{equation}\label{eq3.19}
\int_{S^{n-1}\cap\varkappa(\Omega)}|\nabla_\omega
v(t,\omega)|^2\,d\omega\geq (n-1)
\int_{S^{n-1}\cap\varkappa(\Omega)}|v(t,\omega)|^2\,d\omega,
\end{equation}

\noindent  and the same estimate holds with $v$ replaced by
$\partial_t v$. Together with (\ref{eq3.13}) and (\ref{eq3.14})
this gives

\begin{eqnarray}\label{eq3.20}
&&\hskip -1.7cm \int_{\Omega}\Delta
u(x)\Delta\left(\frac{u(x)g(\log (|\xi|/|x|))}{|x|^{n}}\right)\,dx\nonumber\\[4pt]
&&\hskip -1.7cm\quad \geq\int_{\varkappa(\Omega)}\Bigl[g(t-\tau)
\Bigl( (\delta_\omega v)^2
-2n(\nabla_\omega v)^2 \Bigr)\nonumber\\[4pt]
&&\hskip -1.7cm\qquad \quad
-\Bigl(2\partial_t^2g(t-\tau)+3n\partial_t
g(t-\tau)+\left(n^2-2\right)g(t-\tau)\Bigr)(\partial_t v)^2
\nonumber\\[4pt]
&&\hskip -1.7cm\qquad +\frac
12\Bigl(\partial_t^4g(t-\tau)+2n\partial_t^3g(t-\tau)
+\bigl(n^2-2\bigr)\partial_t^2g(t-\tau)
-2n\partial_tg(t-\tau)\Bigr)v^2\Bigr]\,d\omega dt.
\end{eqnarray}

\noindent On the other hand, for every $t\in\RR$

\begin{equation}\label{eq3.21}
\int_{S^{n-1}\cap\varkappa(\Omega)} \left( (\delta_\omega
v(t,\omega))^2 -2n(\nabla_\omega v(t,\omega))^2\right)\,d\omega\geq
0
\end{equation}

\noindent by Lemma~\ref{l3.2.2}, so it remains to estimate the
terms in the last line of (\ref{eq3.20}). However, according to
Lemma~\ref{l3.2}, we have

\begin{eqnarray}\label{eq3.21.1}
\nonumber &&\hskip
-1.5cm\int_{\varkappa(\Omega)}\Bigl(\partial_t^4g(t-\tau)+2n\partial_t^3g(t-\tau)
+\bigl(n^2-2\bigr)\partial_t^2g(t-\tau)
-2n\partial_tg(t-\tau)\Bigr)v^2\Bigr]\,d\omega dt\\[4pt]
&&=\int_{S^{n-1}\cap\varkappa(\Omega)} v^2(\tau,\omega)\,d\omega,
\end{eqnarray}

\noindent which completes the proof.
  \ep

\section{Global estimates: part II}
\setcounter{equation}{0}

\begin{lemma}\label{l4.1}
Suppose $\Omega$ is a bounded Lipschitz domain in $\RR^n$, $O\in
\RR^n\setminus \Omega$, and

\begin{equation}\label{eq4.1}
u\in C^4(\bar \Omega),\quad u\Bigl|_{\po}=0,\quad\nabla
u\Bigl|_{\po}=0,\quad v=e^{2 t}(u\circ \varkappa^{-1}).
\end{equation}

\noindent Then
\begin{eqnarray}\label{eq4.2}
&&2 \int_{\Omega}\Delta u(x)\Delta\left(\frac{u(x)G(\log
|x|^{-1})}{|x|^{n}}\right)\,dx-\int_{\Omega}\Delta^2
u(x)\,\left(\frac{\big(x\cdot \nabla u(x)\big)\,G(\log
|x|^{-1})}{|x|^{n}}\right)\,dx\nonumber\\[4pt]
&&\nonumber\\[4pt]
&&\quad =\int_{\varkappa(\Omega)}\Bigl(-{\textstyle{\frac 12}}
(\delta_\omega v)^2\partial_tG+ (\partial_t\nabla_\omega
v)^2\left(\partial_t
G+2nG\right)\nonumber\\[4pt]
&&\qquad + (\partial_t^2v)^2\left({\textstyle{\frac 32}}\,
\partial_tG+2nG\right) +n(\nabla_\omega v)^2\partial_t G
\nonumber\\[4pt]
&&\qquad +(\partial_t v)^2\left(-{\textstyle{\frac 12}}
\partial_t^3 G-n\partial_t^2G-{\textstyle{\frac 12}}
 (n^2+2n-4)\partial_t G-2n(n-2)G\right)\Bigr)\,d\omega dt
\nonumber\\[4pt]
&&\quad  -\frac 12\int_{\varkappa(\po)}\left((\delta_\omega
v)^2+2(\partial_t\nabla_\omega v)^2+(\partial_t^2
v)^2\right)\,G\,\cos (\nu,t)\,d\sigma_{\omega,t},
\end{eqnarray}

\noindent where $\nu$ stands for an outward unit normal to
$\Omega$ and $G$ is a function on $\RR$ for which both sides of
{\rm{(\ref{eq4.2})}} are well-defined.
\end{lemma}

\bp Passing to the coordinates $(t,\omega)$ and using
(\ref{eq2.5}), one can see that

\begin{eqnarray}\label{eq4.3}
&&\int_{\Omega}\Delta^2 u(x)\left(\frac{\big(x\cdot \nabla
u(x)\big)\,G(\log
|x|^{-1})}{|x|^{n}}\right)\,dx \nonumber\\[4pt]
&&\qquad = \int_{\varkappa(\Omega)}
\Lambda(\partial_t,\delta_\omega)
\Lambda(\partial_t-2,\delta_\omega)v\,(2 v-\partial_tv)\,G\,d\omega
dt,
\end{eqnarray}

\noindent and therefore,

\begin{eqnarray}\label{eq4.4}
&&2\int_{\Omega}\Delta u(x)\Delta\left(\frac{u(x)G(\log
|x|^{-1})}{|x|^{n}}\right)\,dx-\int_{\Omega}\Delta^2
u(x)\left(\frac{\big(x\cdot \nabla u(x)\big)\,G(\log
|x|^{-1})}{|x|^{n}}\right)\,dx\nonumber\\[4pt]
&&\qquad = \int_{\varkappa(\Omega)}
\Lambda(\partial_t,\delta_\omega)
\Lambda(\partial_t-2,\delta_\omega)v\,\partial_t v \,G\,d\omega dt\nonumber\\[4pt]
&&\qquad =\int_{\varkappa(\Omega)}
\left(\partial_t^2-(n-2)\partial_t+\delta_\omega\right)
\left(\partial_t^2-(n+2)\partial_t+2n+\delta_\omega\right)v\,\partial_t
v \,G\,d\omega dt.
\end{eqnarray}

\noindent This, in turn, is equal to

\begin{eqnarray}\label{eq4.5}
&\nonumber
&\int_{\varkappa(\Omega)}\Bigl(\delta_\omega^2v+2\delta_\omega\partial_t^2v+\partial_t^4
v-2n\partial_t\delta_\omega v-2n \partial_t^3v-2n(n-2)\partial_t v
\\[4pt]
&&\qquad +2n\delta_\omega v+(n^2+2n-4)\partial_t^2v \Bigr)
\,\partial_t v
\,G\,d\omega dt\nonumber \\[4pt]
&&=\int_{\varkappa(\Omega)} \Bigl(\delta_\omega^2v\partial_t v
G-2\partial_t^2\nabla_\omega v\cdot\partial_t\nabla_\omega v
G-\partial_t^3 v\partial_t^2 v G-\partial_t^3 v\partial_t
v\partial_tG\nonumber \\[4pt]
&&\qquad+2n(\partial_t\nabla_\omega v)^2G+2n
(\partial_t^2v)^2G+2n\partial_t^2 v\partial_t v\partial_t
G-2n(n-2)(\partial_t v)^2 G
\nonumber \\[4pt]
&&\qquad -2n\partial_t\nabla_\omega v\cdot\nabla_\omega v
G+(n^2+2n-4)\partial_t^2v\partial_tv G\Bigr)\,d\omega dt .
\end{eqnarray}

Let us consider the first term in (\ref{eq4.5}). Since the first
derivatives of $v$ vanish on the boundary,

\begin{eqnarray}\label{eq4.6}
\nonumber && \int_{\varkappa(\Omega)} \delta_\omega^2v\partial_t v\,
G\,d\omega dt = -\int_{\varkappa(\Omega)} \nabla_\omega
\delta_\omega
v\cdot\nabla_\omega \partial_t v\, G\,d\omega dt\\[4pt]
\nonumber && \quad =\int_{\varkappa(\Omega)}
\left(-\delta_\omega(\nabla_\omega v\cdot \partial_t\nabla_\omega
v)+\nabla_\omega v\cdot \partial_t\delta_\omega\nabla_\omega
v+2\delta_\omega v\partial_t\delta_\omega v\right) G\,d\omega dt
\\[4pt]
&& \quad =\frac 12\int_{\varkappa(\Omega)}
\left(-\delta_\omega\partial_t(\nabla_\omega v)^2+
\partial_t(\delta_\omega v)^2\right)
G\,d\omega dt.
\end{eqnarray}

\noindent Integrating by parts in $t$, we can rewrite the last
expression as

\begin{eqnarray}\label{eq4.6.1}
\nonumber && \hskip -1cm \frac 12\int_{\varkappa(\Omega)}
\left(\delta_\omega(\nabla_\omega v)^2- (\delta_\omega
v)^2\right)\partial_t G\,d\omega dt-\frac 12\int_{\varkappa(\po)}
\left(\delta_\omega(\nabla_\omega v)^2- (\delta_\omega v)^2\right)
G\,\cos(\nu,t)\,d\sigma_{\omega,t}\\[4pt]
\nonumber && \hskip -1cm =\frac 12\int_{\varkappa(\Omega)}
\left(2\delta_\omega\nabla_\omega v\cdot\nabla_\omega v+
(\delta_\omega v)^2\right)\partial_t G\,d\omega dt-\frac
12\int_{\varkappa(\po)} (\delta_\omega v)^2
G\,\cos (\nu,t)\,d\sigma_{\omega,t}\\[4pt]
 && \hskip -1cm =-\frac 12 \int_{\varkappa(\Omega)}
(\delta_\omega v)^2\partial_t G\,d\omega dt-\frac
12\int_{\varkappa(\po)} (\delta_\omega v)^2 G\,\cos
(\nu,t)\,d\sigma_{\omega,t}.
\end{eqnarray}

Now (\ref{eq4.2}) follows directly from (\ref{eq4.5}) by
integration by parts and (\ref{eq4.6})--(\ref{eq4.6.1}). \ep

\begin{lemma}\label{l4.2}
Suppose $\Omega$ is a bounded convex domain in $\RR^n$, $O\in\po$,

\begin{equation}\label{eq4.7} u\in C^4(\bar \Omega),\quad
u\Bigl|_{\po}=0,\quad\nabla u\Bigl|_{\po}=0,\quad v=e^{2 t}(u\circ
\varkappa^{-1}),
\end{equation}

\noindent  and $g$ is defined by {\rm{ (\ref{eq3.9})}}. Then
\begin{eqnarray}\label{eq4.8}
&&\hskip -1cm 2 \int_{\Omega}\Delta u(x)\Delta\left(\frac{u(x)g(\log
|\xi|/|x|)}{|x|^{n}}\right)\,dx-\int_{\Omega}\Delta^2
u(x)\left(\frac{\big(x\cdot \nabla u(x)\big)\,g(\log
|\xi|/|x|)}{|x|^{n}}\right)\,dx\nonumber\\[4pt]
&&\hskip -1cm \quad \geq -\frac 12 \int_{\varkappa(\Omega)}\Bigg(
\partial_t^3g(t-\tau)+2n\,\partial_t^2 g(t-\tau)\nonumber\\[4pt]
&& \qquad\qquad +(n^2-2)\,\partial_tg(t-\tau) -4n
g(t-\tau)\Bigg)(\partial_t v(t,\omega))^2\,d\omega dt,
\end{eqnarray}

\noindent for every $\xi\in\Omega$, $\tau=\log |\xi|^{-1}$.
\end{lemma}

\bp Observe that $g\geq 0$ and for every convex domain $\cos
(\nu,t)\leq 0$ so that the boundary integral in (\ref{eq4.2}) is
non-positive.

Going further, the formula (\ref{eq3.15}) shows that $\partial _t
g\leq 0$. For $t<0$ it is obvious and when $t>0$ this function can
change sign at most once, while at $0$ and in the neighborhood of
$+\infty$ it is negative. In combination with Lemma~\ref{l3.2.2}
this yields

\begin{equation}\label{eq4.9}
\int_{\varkappa(\Omega)}\left(-{\textstyle{\frac 12}} (\delta_\omega
v(t,\omega))^2\partial_tg(t-\tau) +n(\nabla_\omega
v(t,\omega))^2\partial_t g(t-\tau) \right)\,d\omega dt\geq 0.
\end{equation}

Using (\ref{eq3.9}), (\ref{eq3.15}) we can also show that the
coefficients of $(\partial_t^2 v)^2$ and $(\partial_t\nabla_\omega
v)^2$ are nonnegative. Indeed, we compute

\begin{eqnarray*}
&&\frac 32\,\partial_t\,g(t)+2n\,g(t)\nonumber\\[4pt]
&&\qquad =\frac 1{2\sqrt{n^2+8}}\left\{\begin{array}{l}
\left(-\frac 54 n-\frac 34\sqrt{n^2+8}\right) e^{- 1/2(n-\sqrt{n^2+8})t}+2\sqrt{n^2+8},\qquad\qquad t<0,\\[4pt]
\frac 12\sqrt{n^2+8}\, e^{-nt}+ \left(-\frac 54 n+\frac 34
\sqrt{n^2+8}\right)e^{- 1/2(n+\sqrt{n^2+8})t},\qquad t>0.
\end{array}
\right.
\end{eqnarray*}

\noindent The positivity of this function can be proved in a way
similar to the argument for $\partial_tg$: this time, the function
is positive at $0$ and in the neighborhood of $\pm \infty$. Since
$\partial_t g\leq 0$, we also have that $\partial_tg+2n g\geq 0$.

Finally, since $\Omega$ is a convex domain, we can apply
(\ref{eq3.19}) for the function $\partial_t v(t,\cdot)$,
$t\in\RR$, and obtain the estimate (\ref{eq4.8}).\ep

\begin{lemma}\label{l4.3}
Suppose $\Omega$ is a bounded convex domain in $\RR^n$, $O\in\po$,

\begin{equation}\label{eq4.10}
u\in C^4(\bar \Omega),\quad u\Bigl|_{\po}=0,\quad\nabla
u\Bigl|_{\po}=0,\quad v=e^{\lambda t}(u\circ \varkappa^{-1}).
\end{equation}

\noindent and $g$ is given by {\rm{(\ref{eq3.9})}}. Then

\begin{eqnarray}\label{eq4.11}
\frac 12\int_{S^{n-1}\cap
\varkappa(\Omega)}v^2(\tau,\omega)\,d\omega&\leq & \frac{n^2+n-2}{n}
\int_{\Omega}\Delta
u(x)\Delta\left(\frac{u(x)g(\log |\xi|/|x|)}{|x|^{n}}\right)\,dx \nonumber\\[4pt]
&&-\frac {n^2-2}{2n}\int_{\Omega}\Delta^2
u(x)\left(\frac{\big(x\cdot \nabla
u(x)\big)\,g(\log |\xi|/|x|)}{|x|^{n}}\right)\,dx\nonumber\\[4pt]
\end{eqnarray}

\noindent for every $\xi\in\Omega$, $\tau=\log |\xi|^{-1}$.
\end{lemma}

\bp First of all, (\ref{eq3.12}) combined with (\ref{eq4.8})
implies

\begin{eqnarray}\label{eq4.12}
&&\hskip -0.7 cm (n^2-2)\left(2\int_{\Omega}\Delta u(x)\Delta
\left(\frac{u(x)g(\log |\xi|/|x|)}{|x|^{n}}\right)\,dx
-\int_{\Omega}\Delta^2 u(x)\frac{\big(x\cdot \nabla
u(x)\big)\,g(\log |\xi|/|x|)}{|x|^{n}}\,dx\right)\nonumber\\[4pt]
&& \hskip -0.7 cm\qquad\qquad +2n \int_{\Omega}\Delta u(x)\Delta
\left(\frac{u(x)g(\log |\xi|/|x|)}{|x|^{n}}\right)\,dx
\nonumber\\[4pt]
&& \hskip -0.7 cm\quad\geq n\int_{S^{n-1}\cap \varkappa(\Omega)}v^2(\tau,\omega)\,d\omega
-\int_{\varkappa(\Omega)}\Bigl(\left(n^2/2-1\right)\partial_t^3g(t-\tau)\nonumber\\[4pt]
&& \hskip -0.7 cm\qquad\qquad
+(n^3+2n)\partial_t^2g(t-\tau)+\left(n^4/2+4n^2+2\right)\partial_tg(t-\tau)\Bigr)(\partial_t
v)^2\,d\omega dt\end{eqnarray}

\noindent Using the formulas (\ref{eq3.9}), (\ref{eq3.15}),
(\ref{eq3.16}) and

\begin{equation}\label{eq4.12.1}
\partial_t^3\,g(t)=\frac 1{2\sqrt{n^2+8}}\left\{\begin{array}{l}
 \frac{1}{8}\left(n-\sqrt{n^2+8}\right)^3 e^{-
1/2(n-\sqrt{n^2+8})t},\qquad
\,\,\,\,\,\qquad \qquad\qquad t<0,\\[4pt]
-n^2 \sqrt{n^2+8}\, e^{-nt}+
\frac{1}{8}\left(n+\sqrt{n^2+8}\right)^3e^{-
1/2(n+\sqrt{n^2+8})t},\,\, t>0.
\end{array}
\right.
\end{equation}

\noindent we compute

\begin{eqnarray}
&&\nonumber\hskip -0.7 cm\left(n^2/2-1\right)\partial_t^3g(t)
+(n^3+2n)\partial_t^2g(t)+\left(n^4/2+4n^2+2\right)\partial_tg(t)
=-\frac{1}{2\sqrt{n^2+8}}\times\nonumber\\[4pt]
&& \hskip -0.7 cm\quad
\times\left\{\begin{array}{l}n\left(n\sqrt{n^2+8}+6\right)
e^{-1/2(n-\sqrt{n^2+8})t},\,\,\qquad\qquad\qquad\qquad\qquad\quad t<0,\\[4pt]
(2+n^2)\sqrt{n^2+8} \,e^{-nt}
+n\left(-n\sqrt{n^2+8}+6\right)e^{-1/2(n+\sqrt{n^2+8})t},\,\, t>0.
\end{array}
\right.\label{eq4.12.2}
\end{eqnarray}

\noindent The function above is non-positive,  which in
combination with (\ref{eq4.12}) implies (\ref{eq4.11}). \ep

\section{Local estimates for the solution of the Dirichlet \break problem}
\setcounter{equation}{0}

Throughout this section we will adopt the following notation:

\begin{eqnarray*}
S_r(Q)&:=& \{x\in\RR^n:\,|x-Q|=r\},\qquad S_r:=S_r(O),\\[4pt]
B_r(Q)&:=& \{x\in\RR^n:\,|x-Q|<r\},\qquad B_r:=B_r(O),\\[4pt]
C_{r,R}(Q)&:=& \{x\in\RR^n:\,r\leq |x-Q|\leq R\},\,
C_{r,R}:=C_{r,R}(O),
\end{eqnarray*}

\noindent where $Q\in\RR^n$ and $0<r<R<\infty$.

Having this at hand, let us start with a suitable version of the
Caccioppoli inequality for the biharmonic equation.

\begin{lemma}\label{l5.1} Let $\Omega$ be an arbitrary
domain on $\RR^n$, $Q\in\po$ and $R\in(0,{\rm diam}\,(\Omega)/5)$.
Suppose

\begin{equation}\label{eq5.1}
\Delta^2 u=f \,\,{\mbox{in}}\,\,\Omega, \quad f\in
C_0^{\infty}(\Omega\setminus B_{5R}(Q)),\quad u\in \ring
W_2^2(\Omega).
\end{equation}

\noindent Then

\begin{equation}\label{eq5.2}
\int_{B_{\rho}(Q)}|\nabla^2 u|^2\,dx
+\frac{1}{\rho^2}\int_{B_{\rho}(Q)}|\nabla u|^2\,dx\leq
\frac{C}{\rho^4}\int_{C_{\rho,2\rho}(Q)}|u|^2\,dx
\end{equation}

\noindent for every $\rho<4R$.
\end{lemma}

We now proceed with the local estimates for solutions near a
boundary point of the domain.

\begin{theorem}\label{l5.2}
Let $\Omega$ be a bounded smooth convex domain in $\RR^n$,
$Q\in\po$, and $R\in(0,{\rm diam}\,(\Omega)/5)$. Suppose

\begin{equation}\label{eq5.3}
\Delta^2 u=f \,\,{\mbox{in}}\,\,\Omega, \quad f\in
C_0^{\infty}(\Omega\setminus B_{5R}(Q)),\quad u\in \ring
W_2^2(\Omega).
\end{equation}

\noindent Then

\begin{equation}\label{eq5.4}
\frac{1}{\rho^{4}}\meanint_{S_{\rho}(Q)}|u(x)|^2\,d\sigma_x \leq
\frac{C}{R^4} \meanint_{C_{R,4R}(Q)} |u(x)|^2\,dx\quad {\mbox{ for
every}}\quad \rho<R,
\end{equation}

\noindent where the constant $C$ depends on the dimension only.
\end{theorem}

\bp Assume for the moment that $Q=O$. For a smooth domain $\Omega$
the solution to the boundary value problem (\ref{eq5.3}) belongs
to the class $C^4(\bar\Omega)$ by the standard elliptic theory.

Going further,  take some $\eta\in C_0^\infty(\RR)$ such that
$0\leq\eta\leq 1$ and

\begin{equation}\label{eq5.5}
\eta=0\,\, {\rm for}\,\,t\leq \log(2R)^{-1}, \quad \eta=1\,\, {\rm
for}\,\,t\geq \log R^{-1},\quad |\partial_t^k \eta|\leq C, \quad
k\leq 4.
\end{equation}

Since $u\in C^4(\bar\Omega)$ we can apply Lemma~\ref{l4.3} with
$(\eta\circ\varkappa) u$ in place of $u$. The function $u$ is
biharmonic on the support of $\eta\circ\varkappa$, so that

\begin{equation}\label{eq5.6}
\left[\Delta^2,\eta\circ\varkappa\right]u=\Delta^2((\eta\circ\varkappa)\,
u)-(\eta\circ\varkappa)\Delta^2 u=\Delta^2((\eta\circ\varkappa) u)
\end{equation}

\noindent and hence by (\ref{eq4.11}),

\begin{eqnarray}\label{eq5.7}
&&\hskip -1cm \frac 12 \,\eta^2(\tau)\int_{S^{n-1}\cap \varkappa(\Omega)}v^2(\tau,\omega)\,d\omega\nonumber\\[4pt]
&&\hskip -1cm\quad\leq
\frac{n^2+n-2}{n}\int_{\Omega}\left[\Delta^2,\eta(\log
|x|^{-1})\right] u(x)\left(\frac{u(x)\eta(\log |x|^{-1})g(\log
(|\xi|/|x|))}{|x|^{n}}\right)\,dx\nonumber\\[4pt]
&&\hskip -1cm\quad-\frac
{n^2-2}{2n}\int_{\Omega}\left[\Delta^2,\eta(\log |x|^{-1})\right]
u(x)\left(\frac{x\cdot \nabla (u(x)\eta(\log |x|^{-1}))\,g(\log
(|\xi|/|x|))}{|x|^{n}}\right)\,dx
\nonumber\\[4pt]
&&\hskip -1cm\quad\leq \sum_{i=1}^2\sum_{j,k=0}^2
C_{i,j,k}\int_{\varkappa(\Omega)}(\partial_t^k\nabla_{\omega}^j
v)^2\left(\partial_t^i\eta\right)^2 \,d\omega dt,
\end{eqnarray}

\noindent where $C_{i,j,k}$ are some constants depending on the
dimension only. Here for the last inequality we used integration
by parts, Cauchy-Schwartz inequality and boundedness of the
function $g$ and its derivatives. Observe that
$\left[\Delta^2,\eta\circ\varkappa\right]u$ contains only the
derivatives of $u$ of order less than or equal to 3, for that
reason there are no boundary terms coming from the integration by
parts and the order of derivatives in the final expression does
not exceed 2.

Observe also that  each term on the right hand side of
(\ref{eq5.7}) contains some derivative of $\eta$. However,

\begin{equation}\label{eq5.8}
{\rm supp}\,\partial_t^k\eta\subset\left(\log (2R)^{-1},\,\log
R^{-1}\right), \qquad \mbox{for}\quad k\geq 1,
\end{equation}

\noindent so that (\ref{eq5.7}) entails the estimate

\begin{equation}\label{eq5.10}
\frac{1}{\rho^{4}}\meanint_{S_{\rho}}|u(x)|^2\,d\sigma_x \leq C
\sum_{k=0}^2 \frac{1}{R^{4-2k}}\meanint_{C_{R,2R}} |\nabla^k
u(x)|^2\,dx,
\end{equation}

\noindent for every $\rho<R$. Next, invoking Lemma~\ref{l5.1}, we
deduce the estimate (\ref{eq5.4}) for $Q=O$. Then (\ref{eq5.4})
follows from it in full generality since $C$ is a constant
depending solely on the dimension. \ep

Given Theorem~\ref{l5.2} it is a matter of approximation to prove
the following result.

\begin{corollary}\label{l5.3}
Let $\Omega$ be a bounded convex domain in $\RR^n$, $Q\in\po$, and
$R\in(0,{\rm diam}\,(\Omega)/10)$. Suppose

\begin{equation}\label{eq5.11}
\Delta^2 u=f \,\,{\mbox{in}}\,\,\Omega, \quad f\in
C_0^{\infty}(\Omega\setminus B_{10R}(Q)),\quad u\in \ring
W_2^2(\Omega).
\end{equation}

\noindent Then

\begin{equation}\label{eq5.12}
\frac{1}{\rho^{4}}\meanint_{C_{\rho/2,\rho}(Q)}|u(x)|^2\,dx \leq
\frac{C}{R^4} \meanint_{C_{R/2,5R}(Q)} |u(x)|^2\,dx\quad {\mbox{ for
every}}\quad \rho<R/2,
\end{equation}

\noindent where the constant $C$ depends on the dimension only.
\end{corollary}

\bp Let us start approximating $\Omega$ by a sequence of smooth
convex domains $\{\Omega_n\}_{n=1}^\infty$ such that

\begin{equation}\label{eq5.13}
\bigcup_{n=1}^\infty \Omega_n=\Omega,\quad
{\overline{\Omega}}_n\subset\Omega_{n+1} \quad \mbox{for every}
\quad n\in\NN.
\end{equation}

\noindent  Choose $N_0\in\NN$ such that ${\rm supp}\,f\subset
\Omega_n$ for every $n\geq N_0$ and denote by $u_n$, the solution
of the Dirichlet problem

\begin{equation}\label{eq5.14}
\Delta^2 u_n = f\quad {\rm in} \quad \Omega_n, \quad f\in
C_0^\infty(\Omega_n),\quad u_n\in \ring W_2^2(\Omega_n), \quad
n\geq N_0.
\end{equation}

\noindent Finally, let $Q_n\in\po_n$ be a sequence of points
converging to $Q\in\po$.

Now fix some $\rho$ as in (\ref{eq5.12}) and let
$\alpha=2^{-100}$. Then there exists $N=N(\rho)\geq N_0$ such that
$|Q-Q_n|<\alpha\rho$ for every $n>N$. In particular,

\begin{equation}\label{eq5.15}
C_{\rho/2,\,\rho}(Q)\subset
C_{(1/2-\alpha)\rho,\,(1+\alpha)\rho}(Q_n)\quad{\mbox{and}}\quad
C_{R,4R}(Q_n)\subset C_{R-\alpha\rho,\,4R+\alpha\rho}(Q).
\end{equation}

\noindent Therefore, for every $\rho$ such that $(1+\alpha)\rho<R$

\begin{eqnarray}\nonumber
&& \frac{1}{\rho^{2}}\left(\meanint_{C_{\rho/2,\,
\rho}(Q)}|u(x)|^2\,dx\right)^{1/2}\\[4pt]
\nonumber &&\quad  \leq \frac{1}{\rho^{2}}
\left(\meanint_{C_{\rho/2,\,\rho}(Q)}|u(x)-u_n(x)|^2\,dx\right)^{1/2}
+\frac{1}{\rho^{2}}
\left(\meanint_{C_{(1/2-\alpha)\rho,\,(1+\alpha)\rho}(Q_n)}
|u_n(x)|^2\,dx\right)^{1/2}\\[4pt]
&&\quad  \leq \frac{1}{\rho^{2}}
\left(\meanint_{C_{\rho/2,\rho}(Q)}|u(x)-u_n(x)|^2\,dx\right)^{1/2}
+\frac{C}{R^{2}} \left(\meanint_{C_{R,4R}(Q_n)}
|u_n(x)|^2\,dx\right)^{1/2}, \label{eq5.16}
\end{eqnarray}

\noindent where we have used (\ref{eq5.15}) and
Theorem~\ref{l5.2}. Similarly, the last term above can be further
estimated by

\begin{equation}
\frac{C}{R^{2}} \left(\meanint_{C_{R,4R}(Q_n)}
|u_n(x)-u(x)|^2\,dx\right)^{1/2}+\frac{C}{R^{2}}
\left(\meanint_{C_{R-\alpha\rho,\,4R+\alpha\rho}(Q)}
|u(x)|^2\,dx\right)^{1/2}. \label{eq5.17}
\end{equation}

Now recall that the solutions $u_n$, being extended by zero
outside of $\Omega_n$ and treated as elements of $\ring
W_2^2(\Omega)$, strongly converge to $u$ in $\ring W_2^2(\Omega)$
(see, e.g., \cite{Necas}). Then passing to the limit as $n\to
+\infty$ we conclude that

\begin{equation}
\frac{1}{\rho^{2}}\left(\meanint_{C_{\rho/2,
\rho}(Q)}|u(x)|^2\,dx\right)^{1/2} \leq \frac{C}{R^{2}}
\left(\meanint_{C_{R-\alpha\rho,\,4R+\alpha\rho}(Q)}
|u(x)|^2\,dx\right)^{1/2}, \label{eq5.18}
\end{equation}

\noindent for every $\rho<R/(1+\alpha)$ and $\alpha$ sufficiently
small, for instance, $\alpha=2^{-100}$, and finish the argument.
\ep

Finally, we are ready for the \vskip 0.08in

\noindent{\it Proof of Theorem~\ref{t1.1}.} The interior estimates
for solutions of elliptic equations (see \cite{ADN}) imply that
for $x\in B_{R/5}$

\begin{equation}\label{eq5.19}
|\nabla^2 u(x)|^2\leq C \meanint_{B_{d(x)/2}(x) } |\nabla^2
u(y)|^2\,dy
\end{equation}

\noindent where $d(x)$ denotes the distance from $x$ to $\po$.
Denote by $x_0$ a point on the boundary of $\Omega$ such that
$d(x)=|x-x_0|$. Since $x\in B_{R/5}=B_{R/5}(O)$, we have $x\in
B_{R/5}(x_0)$, and therefore for $\alpha=2^{-100}$ we obtain

\begin{displaymath}
\meanint_{B_{d(x)/2}(x)} |\nabla^2 u(y)|^2\,dy\leq
\frac{C}{d(x)^4} \meanint_{C_{d(x)/2,3d(x)}(x_0)} |u(y)|^2\,dy
\leq \frac{C}{R^{4}}\meanint_{C_{(1-\alpha)R,(4+\alpha)R}(x_0)}
|u(y)|^2\,dy,
\end{displaymath}

\noindent using Lemma~\ref{l5.1} for the first estimate and
(\ref{eq5.18}) with $Q=x_0$ for the second one. Since
\begin{equation}\label{eq5.20}
x\in B_{R/5}(O),\quad x\in B_{R/5}(x_0), \quad y\in
C_{(1-\alpha)R,(4+\alpha)R}(x_0) \quad \Longrightarrow \quad y\in
C_{R/2,\,5R}(O),
\end{equation}

\noindent combining (\ref{eq5.19})-- (\ref{eq5.20}) we finish the
proof. \ep

\vskip 0.08in \noindent --------------------------------------
\vskip 0.10in

\noindent {\it Svitlana Mayboroda}

\noindent Department of Mathematics, The Ohio State University,\\
231 W 18th Av., Columbus, OH, 43210, USA\\
{\tt svitlana\@@math.ohio-state.edu}

\vskip 0.10in

\noindent {\it Vladimir Maz'ya}

\noindent Department of Mathematics, The Ohio State University,\\
231 W 18th Av., Columbus, OH, 43210, USA

\vskip 0.05in

\noindent Department of Mathematical Sciences, M\&O Building,\\
University of Liverpool, Liverpool L69 3BX, UK

\vskip 0.05in

\noindent Department of Mathematics, Link\"oping University,\\
SE-581 83 Link\"oping, Sweden

\noindent {\tt vlmaz\@@math.ohio-state.edu,  vlmaz\@@mai.liu.se}

\end{document}